# TEACHING RESOURCES AND TEACHERS' PROFESSIONAL DEVELOPMENT:

## TOWARDS A DOCUMENTATIONAL APPROACH OF DIDACTICS

Ghislaine Gueudet / Luc Trouche

CREAD, IUFM Bretagne UBO / EducTice, INRP; LEPS, Université de Lyon

*In this paper we propose a theoretical approach of teachers' professional development, focusing on teachers' interactions with resources, digital resources in particular. Documents, entailing resources and schemes of utilization of these resources, are developed through documentational geneses occurring along teachers' documentation work (selecting resources, adapting, combining, refining them). The study of teachers' documentation systems permits to seize the changes brought by digital resources; it also constitutes a way to capture teachers' professional change.*

## INTRODUCTION

We present in this paper the first elements of a theoretical approach elaborated for the study of teachers' development, and in particular teachers ICT integration.

Technology integration, and the way teachers work in technology-rich environments, have been extensively researched, and discussed at previous CERME conferences (Drijvers *et al.*, 2005, Kynigos *et al.*, 2007). Ruthven's presentation at CERME 5 drew attention on the structuring context of the classroom practice, and on its five key features: *working environment, resource system, activity format, curriculum script, time economy* (Ruthven, 2007). This leads in particular to consider ICT as part of a wider range of available teaching resources. This view also fits technological evolutions: most of paper material is now at some point in digital format; teachers exchange digital files by e-mail, use digital textbooks, draw on resources found on websites etc. Considering ICT amongst other resources raises the question of connections between research on ICT and resources-oriented research.

Many research works address the use of *curriculum material* (Ball & Cohen, 1996; Remillard, 2005). They observe the influence of such material on the enacted curriculum, but also highlight the way teachers shape the material they draw on, introducing a vision of "curriculum use as participation with the text" (Remillard, 2005, p.121). Other authors consider more general resources involved in teaching: material and human, but also mathematical, cultural and social resources (Adler, 2000). They analyze the way teachers interpret and use the available resources, and the consequences of these processes on teachers' professional evolution.

Such statements sound familiar for researchers interested in ICT, who "consider not only the ways in which digital technologies shape mathematical learning through novel infrastructures, but also how it is shaped by its incorporation into mathematical learning and teaching contexts" (Hoyles & Noss, 2008, p. 89). Conceptualization of these processes is provided by the *instrumental approach* (Guin *et al.*, 2005) and by the work of Rabardel (1995) grounding it; this theoretical frame has contributed to set many insightful results about the way students learn mathematics with ICT. Further refinements of this theory have led to take into account the role of the teacher and her intervention on students instrumental geneses, introducing the notion of *orchestration* (Trouche, 2004). Considering instrumental geneses for teachers has been proposed in the context of spreadsheets (Haspekian, 2008) and e-exercises bases (Bueno-Ravel & Gueudet, 2007). These refinements can be considered as first steps towards the introduction of concepts coming from the instrumental approach and illuminating the interactions between teachers and ICT.

Thus connections between studies about the use of teaching resources, and studies about the way in which teachers work in a technology-rich environment exist; however, elaborating a theoretical frame encompassing both perspectives requires a specific care. We present here an approach designed for this purpose, and aiming at studying teachers' documentation work: looking for resources, selecting, designing mathematical tasks, planning their order, carrying them out in class, managing the available artefacts, etc. We take into account teachers' work in class, but also their (too often neglected) work out of class.

We draw on the theoretical elements evoked above, but also on field data. Some of these data come from previous research in which we were engaged: particularly about use of e-exercises bases (Bueno-Ravel & Gueudet, 2007), and about an in-service training design, the SFoDEM (Guin & Trouche, 2005). Other data were specifically collected: we have set up a series of interviews with nine secondary school teachers. We chose teachers with different collective involvements, different institutional contexts and responsibilities, and different ICT integration degrees (Assude, 2007). We met them at their homes (where, in France and for secondary teachers, most of their documentation work takes place), and asked them about their uses of resources, and the evolution of these ways of use. We observed the organization of their workplaces at home, of their files (both paper and digital), and collected materials they designed or used. The analyses of these data contributed to shape the concepts; in this paper we only use them to display illustrations of the theory. All the interviews took place in France; thus the national context certainly influences the results we display. We hypothesize nevertheless that the concepts exposed are likely to illuminate documentation work in diverse situations.

We present here the elementary concepts of this theory, introducing in particular a distinction between *resources* and *documents*, and the notion of *documentational genesis*. We also expose the specific view of professional evolutions it entails.

## RESOURCES, DOCUMENTS, DOCUMENTATIONAL GENESES

The instrumental approach (Rabardel, 1995, Guin *et al.*, 2005) proposes a distinction between *artefact* and *instrument*. An artefact is a cultural and social means provided by human activity, offered to mediate another human activity. An instrument comes from a process, named *instrumental genesis*, along which the subject builds a *scheme* of utilization of the artefact, for a given class of situations. A scheme, as Vergnaud (1998) defined it from Piaget, is an *invariant organization of activity* for a given class of situations, comprising in particular rules of action, and structured by *operational invariants*, which consist of implicit knowledge built through various contexts of utilization of the artefact. Instrumental geneses have a dual nature. On the one hand, the subject guides the way the artefact is used and, in a sense, *shapes* the artefact: this process is called *instrumentalization*. On the other hand, the affordances and constraints of the artefact influence the subject's activity: this process is called *instrumentation*. We propose here a theoretical approach of teaching resources, inspired by this instrumental approach, with distinctive features that we detail hereafter, and a specific vocabulary.

We use the term *resources* to emphasize the variety of the artefacts we consider: a textbook, software, a student's sheet, a discussion with a colleague etc. A resource is never isolated: it belongs to a set of resources. The subjects we study are teachers. A teacher draws on resources sets for her documentation work. A genesis process takes place, bearing what we call *a document*. The teacher builds schemes of utilization of a set of resources, for the same class of situations, across a variety of contexts. The formula we retain here is:

*Document = Resources + Scheme of Utilization.*

A document entails, in particular, operational invariants, which consist of implicit knowledge built through various contexts of utilization of the artefact, and can be inferred from the observation of invariant behaviors of the teacher for the same class of situations across different contexts.

Figure 1 represents a *documentational genesis*. The instrumentalization process conceptualizes teacher appropriating and reshaping resources, and the instrumentation process captures the influence, on the teacher's activity, of the resources she draws on.

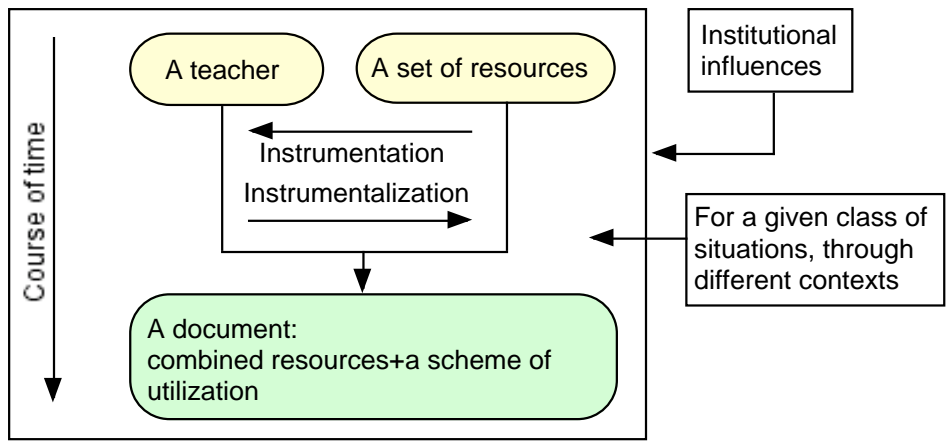

**Figure 1. Schematic representation of a documentational genesis**

## DOCUMENTATIONAL GENESES: TWO ILLUSTRATIVE EXAMPLES

We use a first case study (figure 2) coming from our interviews to illustrate the distinction between a set of resources and a document.

Marie-Pierre (aged 40, involved in collective work within an IREM[1] group; no institutional responsibilities, strong degree of ICT integration) is teaching at secondary school for 14 years, from grade 6 to 9. She uses dynamic geometry systems, spreadsheets, and many online resources (e-exercises and mathematics history websites in particular). She has a digital version of the class textbook. Marie-Pierre has an interactive whiteboard in her classroom for three years, and uses it in each of her courses. For the introduction of the circle's area in grade 7, she starts in class by using a website comprising historical references (Archimedes using circular sections to link the perimeter and the area of a circle) and displaying an animation of the circle unfolding and transforming into a triangle (roughly, but that point is not discussed). Then she presents her own course, based on an extract of the class digital textbook. She complements as usual the files displayed on the whiteboard by writing additional comments and explanations, highlighting important expressions etc.

**Figure 2. Marie-Pierre, example of a lesson introducing the circle's area**

For the class of situations: "design and implement the introduction of the circle's area in grade 7" (figure 2), Marie-Pierre draws on a set of resources comprising the

---

[1] Institute for Research on Mathematics Teaching.

interactive whiteboard, a website[2], a digital textbook, and a hard copy of it. The official curriculum texts, about the circle area, only state that "an inquiry-based approach permits to check the area formula", with no more details. The digital textbook proposes an introductory activity with a digital geometry software: drawing circles, and displaying their areas. Several radius are tested, the radius square and the corresponding area are noted by the students in a table, and they are asked to observe that they obtain an (approximate) ratio table. But Marie-Pierre prefers to draw on the website animated picture (both choices correspond more to an observation activity for the students than to an inquiry-based approach, but we will not discuss this aspect here). So, we claim that she has developed a scheme of utilization of this set of resources, structured by several operational invariants. These invariants are professional beliefs that we infer from our data:

-"A new area formula must be justified by an animation showing a cutting and recombining of the pieces to form a figure whose area is known". This operational invariant concerns all the areas introduced, it also intervenes in the document corresponding to the introduction of the triangle's area for example;

-"The circle's area must be linked with a previously known area: the triangle"; "The circle's area must be linked with the circle's perimeter". These operational invariants are related with the precise mathematical content of the lesson, they were built along the years, with different grade 7 classes (Marie-Pierre uses this website's animation for three years, with two grade 7 classes each year).

We do not assert that these operational invariants were not present among Marie-Pierre's professional knowledge before her integration of the interactive whiteboard. But the possibility to display an animation on a website, to complement it by writing additional explanations, to go back to a previous state of the board to link the "official" formula with what has been observed, yielded a document integrating these operational invariants. And we claim that the development of this document is likely to reinforce, in particular, the above presented beliefs. The operational invariants are both driving forces and outcomes of the teacher's activity.

Documentational geneses are ongoing processes; we use a second case study (figure 3) to emphasize this important aspect. Rabardel & Bourmaud (2005) claim that the design *continues in usage*. We consider here accordingly that a document developed from a set of resources provides new resources, which can be involved in a new set of resources, which will lead to a new document etc. Because of this process, we speak of a *dialectical* relationship between resources and documents.

---

[2] http://pagesperso-orange.fr/therese.eveilleau/pages/hist_mat/textes/mirliton.htm

Marie-Françoise (aged 55; involved in collective work within an IREM group, institutional responsibilities as in-service teacher trainer, strong degree of ICT integration) works with students from grade 10 to 12. She organizes for them 'research narratives': problem solving sessions, where students work in groups on a problem and write down their own 'research narratives' (both solutions and research processes). Thus one class of situations, for Marie-Françoise is 'elaborating open problems for research narratives sessions'. For this class of situations, she draws on a set of resources comprising various websites, but also personal existing resources, colleagues' ideas, etc. Students' ideas constitute a major resource for Marie-Françoise, as she told us: "There is the problem and the way you enact it, because students are free to invent things, and afterwards we benefit from the richness of all these ideas, and you can build on it". The design clearly goes on in class. Moreover, the class sessions provide new resources: the students' research narrative, that Marie-Françoise collects, and saves in a new binder, aiming to enrich the next document built on the same open problem.

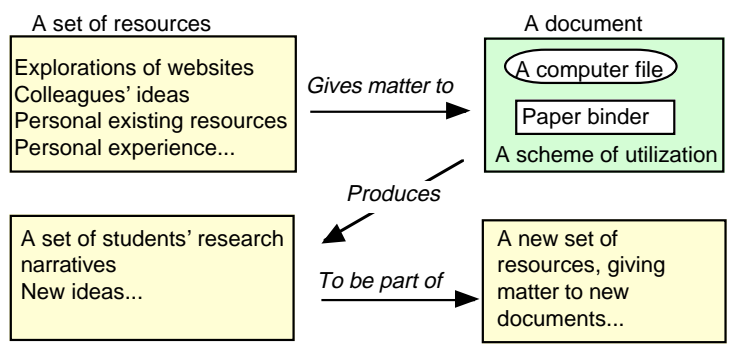

**Figure 3. An illustration of the resources/document dialectical relationship**

The resources evolve, are modified, combined; documents develop along geneses and bear new resources (figure 3) etc. We consider that these processes are part of teachers' professional evolutions, and play a crucial role in them.

## DOCUMENTATION SYSTEMS AND PROFESSIONAL DEVELOPMENT

According to Rabardel (2005), professional activity has a double dimension. Obviously a productive dimension: the outcome of the work done. But the activity also entails a modification of the subject's professional practice and beliefs, within a constructive dimension. Naturally, this modification influences further production processes: the productive/constructive relationship has a dialectical nature.

Teachers' documentation work is the driving force behind documentational geneses, thus it yields productive and constructive professional changes. Literature about teachers' professional change raises the question of the articulation between change of practice and change of knowledge and beliefs. We consider that both are strongly intertwined (e.g., Cooney, 2001). The documentational geneses provide a specific view of this relationship. Working with resources, for the same class of situations across different contexts, leads to the development of a scheme, and in particular of rules of action (professional practice features) and of operational invariants (professional implicit knowledge or beliefs). And naturally these schemes influence the subsequent documentation work. All kinds of professional knowledge are concerned by these processes, the evolutions they generate are not curtailed to curricular knowledge (Schulman, 1986). Thus, studying teachers' documents can be considered as a specific way to study teachers' professional development.

According to Rabardel and Bourmaud (2005), the instruments developed by a subject in his/her professional activity constitute a system, whose structure corresponds to the structure of the subject's professional activity. We hypothesize here similarly that a given teacher develops a structured documentation system.

Let us go back to the example of Marie-Pierre evoked above.

> Marie-Pierre keeps all her "paperboards" (digital files with images corresponding to the successive states of the board). She uses these paperboards at the beginning of a new session, to recall what has been written, by herself or by her students, during the preceding session. On her laptop, Marie-Pierre has one folder for each class level. Each of these folders contains one file with the whole year's schedule, and lessons folders for each mathematical theme. The paperboards are inside the lessons folders. The interactive whiteboard screen below corresponds to the introduction of equations in grade 7, in the context of triangles areas.
>
> 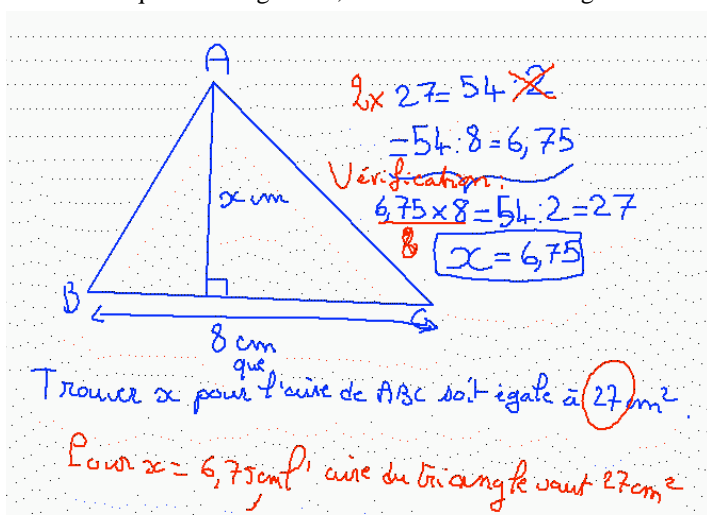
>
> (*Translation:* Find *x* such that ABC area equals 27 cm$^2$. For *x* = 6.75cm, the triangle's area is 27 cm$^2$).

**Figure 4. A view on Marie-Pierre's documents**

Marie-Pierre's files organization on her computer (figure 4), and her statements during the interviews, clearly indicate articulations between her documents. The document whose material component is the year schedule naturally influenced her lesson preparations; but on the opposite, the documents she developed for lessons preparations during previous years certainly intervened in the schedule design. Documents corresponding to connected mathematical themes are also connected. For a given lesson, the students' interventions can contribute to generate operational invariants that will intervene in preparations about other related topics.

A teacher's documents constitute a system, whose organization matches the organization of her professional activity. The evolutions of this documentation system correspond to professional evolutions. Integration of new materials is a visible of the professional practice, and of the documentation system (in the approach we propose, this integration means that a new material is inserted in a set of resources involved in the development of a document). When Marie-Pierre integrates the interactive whiteboard in her courses, it entails a productive dimension: she now teaches with this whiteboard. But it also yields other changes of her practice: she makes more links with previous sessions, in particular recalling students productions is now present in her orchestration choices. And it even generates changes in her professional beliefs, for example about the possible participation of students to her teaching. She seems to have developed an operational invariant like: "a good way to launch a lesson is to recall students' interventions done during the preceding lesson".

The integration of new material is always connected with professional practice and professional beliefs evolutions. But professional evolutions do not always correspond to integration of new material, and the same is true for documentation systems evolutions. For example, Arnaud (47 years old, no collective involvement, institutional responsibilities as in-service teacher trainer, low degree of ICT integration) presented during his interview "help sheets", that he designed years ago for students encountering specific difficulties. He now uses the same sheets as exercises for the whole class; thus while no changes can be observed in the material, the action rules associated evolved.

Integration of new material remains an important issue, especially when the focus is on ICT. The study of a given teacher's documentation system also provides insights in the reasons for the integration or non-integration of a given material. The integration depends indeed on the possibility for this material to be involved in the development of a document, that will be articulated with others within the documentation system. For many years Marie-Pierre prepares her courses as digital files, she uses dynamic geometry software and online resources; the interactive whiteboard articulates with this material. Moreover, Marie-Pierre is convinced of the necessity of fostering students' interventions, and even of including these in the written courses, and the interactive whiteboard matches this conviction. Possible material articulations are important; but other types of articulations must be taken into account, and the integration of new material also strongly depends on operational invariants, thus on teachers' professional knowledge and beliefs.

## CONCLUSION

This paper is related with the second theme of WG7: *Interaction between resources and teachers' professional practice*. It introduces a conceptualization of teachers' interactions with resources and of the associated professional development. Here we just presented the first concepts of a theory whose elaboration is still in progress. Studying teachers' documentation work requires to set specific methodologies, permitting to capture their work in and out of class, to precise their professional beliefs, and to follow long-term processes: it is the main goal of our research. We did not discuss here the very important issue of collective documentation work, which causes particular processes: its study raises the question of collective documentational genesis and documentation systems, and raises new theoretical needs. The documentational approach we propose also needs to be confronted with other teaching contexts: primary school, tertiary level; diverse countries; and also outside the field of mathematics. Further research is clearly needed; the present evolutions of digital resources make it a major challenge for the studies of teachers' professional evolutions.

**REFERENCES**

Adler, J. (2000). Conceptualising resources as a theme for teacher education. *Journal of Mathematics Teacher Education, 3*, 205-224.

Assude, T. (2007). Teachers' practices and degree of ICT integration. In D. Pitta-Pantazi & G. Philippou (Eds.), *Proceedings of the fifth congress of the European Society for Research in Mathematics Education*, CERME 5, Larnaca, Cyprus, http://ermeweb.free.fr/CERME5b/.

Ball, D.L. & Cohen, D. (1996). Reform by the book: what is -or might be- the role of curriculum materials in teacher learning and instructional reform? *Educational researcher, 25*(9), 6-8, 14.

Bueno-Ravel, L. & Gueudet, G. (2007). Online resources in mathematics: teachers' genesis of use. In D. Pitta-Pantazi & G. Philippou (Eds.), *Proceedings of the fifth congress of the European Society for Research in Mathematics Education*, CERME 5, Larnaca, Cyprus, http://ermeweb.free.fr/CERME5b/.

Cooney, T.J. (2001). Considering the paradoxes, perils and purposes for conceptualizing teacher development. In F.L. Lin & T.J. Cooney (Eds.), *Making sense of mathematics teachers education* (pp. 9-31). Dordrecht: Kluwer Academic Publishers.

Drijvers, P., Barzel, B., Maschietto, M. & Trouche, L. (2006). Tools and technologies in mathematical didactics. In M. Bosch (Ed.), *Proceedings of the Fourth Congress of the European Society for Research in Mathematics Education*, CERME 4, San Feliu de Guíxols, Spain, http://ermeweb.free.fr/CERME4/.

Gueudet, G. & Trouche, L. (online). Towards new documentation systems for mathematics teachers? *Educational Studies in Mathematics*, DOI 10.1007/s10649-008-9159-8.

Guin, D., Ruthven, K. & Trouche, L. (Eds.). (2005). *The didactical challenge of symbolic calculators: turning a computational device into a mathematical instrument*. New York: Springer.

Guin, D. & Trouche, L. (2005). Distance training, a key mode to support teachers in the integration of ICT? Towards collaborative conception of living pedagogical resources. In M. Bosch (Ed.), *Proceedings of the Fourth Conference of the European Society for Research in Mathematics Education*, CERME 4, San Feliu de Guíxols, Spain, http://ermeweb.free.fr/CERME4/.

Haspekian, M. (2008). Une genèse des pratiques enseignantes en environnement instrumenté. In F. Vandebrouck (Ed.), *La classe de mathématiques : activité des élèves et pratiques des enseignants* (pp. 293-318). Toulouse: Octarès.

Kynigos, C., Bardini, C., Barzel, B. & Maschietto, M. (2007). Tools and technologies in mathematical didactics. In D. Pitta-Pantazi & G. Philippou (Eds.), *Proceedings of the fifth* *congress of the European Society for Research in Mathematics Education*, CERME 5, Larnaca, Cyprus, http://ermeweb.free.fr/CERME5b/.

Hoyles, C. & Noss, R. (2008). Next steps in implementing Kaput's research programme. *Educational Studies in Mathematics, 68*(2), 85-97.

Rabardel, P. (1995). *Les hommes et les technologies, approche cognitive des instruments contemporains.* Paris: Armand Colin (English version at http://ergoserv.psy.univ-paris8.fr/Site/default.asp?Act_group=1).

Rabardel, P. (2005). Instrument subjectif et développement du pouvoir d'agir *(Subjective instrument and development of action might)*. In P. Rabardel & P. Pastré (Eds.)*, Modèles du sujet pour la conception. Dialectiques activités développement* (pp. 11–29). Toulouse: Octarès.

Rabardel, P. & Bourmaud, G. (2005). Instruments et systèmes d'instruments *(Instruments and systems of instruments)*. In P. Rabardel & P. Pastré (Eds.), *Modèles du sujet pour la conception. Dialectiques activités développement* (pp. 211-229). Toulouse: Octarès.

Remillard, J.T. (2005). Examining key concepts in research on teachers' use of mathematics curricula. *Review of Educational Research, 75*(2), 211-246.

Ruthven, K. (2007). Teachers, technologies and the structures of schooling. In D. Pitta-Pantazi & G. Philippou (Eds.), *Proceedings of the fifth congress of the European Society for Research in Mathematics Education*, CERME 5, Larnaca, Cyprus, http://ermeweb.free.fr/CERME5b/.

Schulman, L. (1986). Those who understand: knowledge growth in teaching. *Educational Researcher, 15*, 4-14.

Trouche, L. (2004). Managing the complexity of human/machine interactions in computerized learning environments: Guiding students' command process through instrumental orchestrations. *International Journal of Computers for Mathematical Learning, 9*, 281-307.

Vergnaud, G. (1998). Toward a cognitive theory of practice. In A. Sierpinska & J. Kilpatrick (Eds.), *Mathematics education as a research domain: a search for identity* (pp. 227-241). Dordrecht: Kluwer Academic Publisher.